\documentclass[twoside, a4paper, 11pt]{article}
\title{Pl{\"u}nnecke inequalities for measure graphs with applications}
\author{Kamil Bulinski and Alexander Fish}
\usepackage{amsfonts}
\usepackage{amsthm}
\usepackage{verbatim}
\usepackage{amsmath, amssymb}
\usepackage{tikz}
\usetikzlibrary{matrix, arrows}
\usepackage{listings}
\usepackage{color}
\usepackage{listings}
\usepackage[all]{xy}
\usepackage[pdftex,colorlinks,linkcolor=blue,citecolor=blue]{hyperref}
\usepackage{graphicx}
\usepackage{float}
\usepackage[margin=3cm]{geometry}
\begin{document}
\maketitle
\raggedbottom

\nocite{TaoVu, Plunnecke1970, BjorklundFishPlunnecke, countableamenablegroups2010, JinPlunnecke, PetridisGraph, Ruzsasumsetandstructure, FurstenbergSzemeredi}

\newcounter{dummy} \numberwithin{dummy}{section}

\theoremstyle{definition}
\newtheorem{mydef}[dummy]{Definition}
\newtheorem{prop}[dummy]{Proposition}
\newtheorem{corol}[dummy]{Corollary}
\newtheorem{thm}[dummy]{Theorem}
\newtheorem{lemma}[dummy]{Lemma}
\newtheorem{eg}[dummy]{Example}
\newtheorem{notation}[dummy]{Notation}
\newtheorem{remark}[dummy]{Remark}
\newtheorem{claim}[dummy]{Claim}
\newtheorem{Exercise}[dummy]{Exercise}

\begin{abstract} We generalize Petridis's new proof of Pl{\"u}nnecke's graph inequality \cite{PetridisGraph} to graphs whose vertex set is a measure space. Consequently, by a recent work of Bj{\"o}rklund and Fish \cite{BjorklundFishPlunnecke}, this gives new Pl{\"u}nnecke inequalities for measure preserving actions which enable us to deduce, via a Furstenberg correspondence principle, Banach density estimates in countable abelian groups that improve on those given by Jin in \cite{JinPlunnecke}. \end{abstract}

\section{Introduction}

\subsection{Background and summary of results}

Given an abelian group $G$ with subsets $A, B \subset G$, it is of great interest to estimate the size of the \textit{product set} (commonly referred to as the \textit{sumset} when additive notation is employed) defined by $$AB = \{ab | a\in A, b \in B\}.$$ In particular, one is also interested in the sizes of iterated product sets $B^k$, which may be recursively defined by $B^1=B$ and $B^k=B^{k-1}B$ for positive integers $k>1$. General inequalities regarding the cardinalities of these were given by Pl{\"u}nnecke and Ruzsa, a comprehensive treatment of which may be found in \cite{Ruzsasumsetsandstructure}. In particular, it was shown in \cite{Plunnecke1970} that if one defines, for finite sets $A \subset G$ and $B \subset G$, the magnification ratios $$D_k = \min_{B' \subset B, B' \neq \emptyset} \frac{|A^kB'|}{|B'|}$$ then $D_k^{1/k}$ is a decreasing sequence in $k$. For infinite subsets we can no longer use cardinalities. It is natural, instead, to use the notion of an invariant density. This notion can be defined in countable amenable groups - the groups which possess a F{\o}lner sequence. By a F{\o}lner sequence, we mean a sequence of finite subsets $F_n \subset G$ such that for each $g \in G$ we have $$\lim_{n \to \infty} \frac{|F_n \cap gF_n|}{|F_n|} = 1.$$ A countable group for which a F{\o}lner sequence exists is called amenable. It is well known that all countable abelian groups are amenable. Each F{\o}lner sequence $(F_n)_{n \in \mathbb{N}}$ gives rise to corresponding upper and lower densities for $A \subset G$ given by $$\overline{d}_{(F_n)} (A) = \limsup_{n \to \infty} \frac{|F_n \cap A|}{|F_n|}$$ and $$ \underline{d}_{(F_n)} (A) = \liminf_{n \to \infty} \frac{|F_n \cap A|}{|F_n|}, $$ respectively. For example, the F{\o}lner sequence $F_n=[1,n] \cap \mathbb{Z}$ in the additive group $\mathbb{Z}$ gives rise to the classical upper and lower \textit{asymptotic} densities in the natural numbers, which are of particular interest in Number Theory. Of interest to us is the \textit{upper} (resp. \textit{lower}) \textit{Banach density}, denoted by $d^*(A)$ (resp. $d_*(A)$), which may be defined as the maximum (resp. minimum) of the set $$ \{ \overline{d}_{(F_n)}(A) | (F_n)_{n \in \mathbb{N}} \text{ is F{\o}lner}  \}.$$ The fact that these extrema are attained in any given countable amenable group may be verified by a simple diagonalisation argument. In fact one can show (see Lemma 3.3 of \cite{countableamenablegroups2010}) the stronger assertion that given any F{\o}lner sequence $(F_n)_{n \in \mathbb{N}}$ and $A \subset G$, there exist $t_n \in G$ such that $$d^*(A)=d_{(t_nF_N)} \overline{d}(A),$$ and likewise for lower Banach density. In particular, in the additive group $\mathbb{Z}^d$ one only needs to look at sequences of finite cubes of strictly increasing cardinality. One may now ask whether Banach densities of sumsets satisfy analogous Pl{\"u}nnecke-Ruzsa type inequalities. Indeed, Jin has proven in \cite{JinPlunnecke} such inequalities for the additive semi-group $\mathbb{Z}_{ \geq 0 }$. 

\begin{thm}[Jin] Suppose $A,B \subset \mathbb{Z}_{\geq 0}$, then $$ d^*(A+B) \geq d^*(kA)^{1/k} d^*(B)^{1-\frac{1}{k}} $$ and $$d_*(A + B) \geq d^*(kA)^{1/k} d_{*}(B)^{1-\frac{1}{k}}.$$

\end{thm}

In \cite{BjorklundFishPlunnecke}, Bj\"orklund and Fish introduced a new dynamical approach for Pl\"unnecke type inequalities. As a result, they extended Jin's theorem to any countable abelian group.

\begin{thm}[Bj{\"o}rklund, Fish]\label{thm: j=1 case of Plunnecke} Suppose $G$ is a countable abelian group and $A,B \subset G$. Then we have $$ d^*(AB) \geq d^*(A^k)^{1/k} d^*(B)^{1-\frac{1}{k}} $$ and $$d_*(AB) \geq d^*(A^k)^{1/k} d_{*}(B)^{1-\frac{1}{k}}.$$

\end{thm}

Our work extends \cite{BjorklundFishPlunnecke} and provides analogous lower bounds for $d^*(A^jB)$ and $d_* (A^jB)$. More precisely, we prove following. 

\begin{thm}~\label{thm: Plunnecke for densities}
Suppose that $G$ is a countable abelian group and $A,B \subset G$. Then for integers $0 < j < k$ we have 
$$d^*(A^jB) \geq d^*(A^k)^{j/k} d^{*}(B)^{1-j/k}$$ and $$d_*(A^jB) \geq d^*(A^k)^{j/k} d_{*}(B)^{1-j/k}.$$
\end{thm}

It seems unclear whether or not our results, even for just $G=(\mathbb{Z},+)$, are attainable from an application of Jin's techniques to the inequality $D_k^{1/k} \leq D_j^{1/j}$. We are also able to obtain a result involving multiple different factors, which is a generalization of Theorem~\ref{thm: j=1 case of Plunnecke}. \begin{thm}~\label{thm: Density Plunnecke for different summands} Suppose $G$ is a countable abelian group and $B, A_1, \ldots ,A_k \subset G$. Then $$ d^*(A_1\ldots A_k) \leq d^*(B)^{1-k}\prod_{i=1}^k d^*(A_iB) $$ and $$d^*(A_1\ldots A_k) \leq d_*(B)^{1-k}\prod_{i=1}^k d_*(A_iB) .$$

\end{thm}

The strategy of the proofs of the main theorems involves employing an ergodic approach. This approach was developed by Bj{\"o}rklund and the second author in \cite{BjorklundFishPlunnecke}. First, we prove a Pl{\"u}nnecke inequality for measure preserving actions and then we combine it with a Furstenberg correspondence principle for product sets. Next, we recall the magnification ratios defined for the dynamical setting in \cite{BjorklundFishPlunnecke}.

\begin{mydef}[$G$ acting on a measure space] We say that a group $G$ acts on a measure space $(X, \mathcal{B},\mu)$ if, for each $g \in G$, the map $x \mapsto g.x$ is measure preserving, i.e., it is measurable and $\mu(gB)=\mu(B)$ for each $B \in \mathcal{B}$.

\end{mydef}

\begin{mydef}[\cite{BjorklundFishPlunnecke}] Given a countable abelian group $G$ acting on a measure space $(X, \mathcal{B},\mu)$, let us define, for $A \subset G$ and $B \in \mathcal{B}$ of finite positive measure, the magnification ratio $$c(A,B) = \inf \left \{ \frac{\mu(AB')}{\mu(B')} | B' \subset B, \mu(B')>0 \right \}.$$

\end{mydef}

The following is an extension to general $j<k$ of the result in \cite{BjorklundFishPlunnecke}.

\begin{thm}\label{thm intro: Plunnecke finite translates} If $A$ is finite, then for positive integers $j<k$ we have $$c(A^k,B)^{1/k} \leq c(A^j,B)^{1/j}.$$ 

\end{thm}

Note that the classical Pl\"unnecke-Ruzsa inequality is the case where $X=G$, $\mu$ is the counting measure and the action is by multiplication. It follows from the techniques developed in \cite{BjorklundFishPlunnecke} that Theorem~\ref{thm intro: Plunnecke finite translates} implies Pl\"unnecke inequalities for $A \subset G$ not necessarily finite.

\begin{thm}[Pl{\"u}nnecke inequalities]\label{thm: countable Plunnecke for G system} If $\mu(X)=1$ and $L^2(X,\mathcal{B},\mu)$ is separable, then for positive integers $j<k$ we have $$c(A^k,B)^{1/k} \leq c(A^j,B)^{1/j}.$$

\end{thm}

Next, we recall Furstenberg's correspondence principle for product sets. This correspondence principle can be derived from the seminal work of Furstenberg \cite{FurstenbergSzemeredi}. Nevertheless, it was noticed much later in \cite{countableamenablegroups2010}. The following version of the correspondence principle (the third and fourth inequalities) for product sets is due to Bj\"orklund and the second author \cite{BjorklundFishProductSet}. 

\begin{prop}[Furstenberg's correspondence principle] Suppose that $G$ is a countable abelian group and $A,B \subset G$. Then there exists a compact metrizable space $X$ on which $G$ acts by homeomorphisms such that there exist $G$-invariant ergodic Borel probability measures $\mu, \nu$ on $X$ together with a clopen $\widetilde{B} \subset X$ such that \begin{align*}  d^* (B) &= \mu(\widetilde{B}) \\ d^*(AB) &\geq \mu(A\widetilde{B}) \\ d_* (B)  &\leq \nu(\widetilde{B}) \\ d_*(AB)  &\geq \nu(A\widetilde{B}). \end{align*} 

\end{prop}

Next, we demonstrate how Theorem~\ref{thm: countable Plunnecke for G system}, through Furstenberg's correspondence principle, implies Theorem~\ref{thm: Plunnecke for densities} and Theorem~\ref{thm: Density Plunnecke for different summands}. 

\textbf{Proof of Theorem~\ref{thm: Plunnecke for densities} and Theorem~\ref{thm: Density Plunnecke for different summands}:}  Let $(X,\mu)$ and $\widetilde{B}$ be as in the correspondence principle. Note that we may assume that $d^*(B)>0$ as the result is trivial otherwise. Note also that in Section~\ref{section: Applications to Banach Density estimates} (see Lemma~\ref{lemma: Standard density estimate from ergodic action}) we show that $d^*(A^k) \leq \mu(A^k \widetilde{B})$. Altogether, this gives \begin{align*}\left( \frac{d^*(A^jB)}{d^*(B)} \right)^{1/j} & \geq \left( \frac{\mu(A^j\widetilde{B})}{\mu(\widetilde{B})}  \right)^{1/j} \\ & \geq c(A^j,\widetilde{B})^{1/j}\\ & \geq c(A^k,\widetilde{B})^{1/k} \\ & \geq \left( \frac{d^*(A^k)}{d^*(B)} \right)^{1/k}, \end{align*} which shows the first inequality. The second one may be deduced from the same argument applied to the measure $\nu$ instead of $\mu$ from the correspondence principle. Moreover, Theorem~\ref{thm: Density Plunnecke for different summands} may be obtained from applying the correspondence principle to the Pl{\"u}nnecke inequality for different summands (Proposition~\ref{prop: different summands}). $\blacksquare$

\subsection{Outline of paper}

The main object introduced in this paper is, what we call, a \textit{measure graph}. Section~\ref{section: Definitions} provides all the relevant definitions and basic properties. Intuitively, a measure graph is a directed edge-labelled graph equipped with a measure on the vertex set that mimics certain elementary combinatorial properties of the classical graph-theoretic notion of a \textit{matching}. The aim is to prove a measure-theoretic version of the classical Pl{\"u}nnecke inequality for commutative graphs. The classical approach employs Menger's theorem, which has no obvious measure theoretic analogue. However, Petridis \cite{PetridisGraph} has recently found a new proof of this inequality that avoids the use of Menger's theorem. In Section~\ref{section: Plunnecke's inequality for measure graphs} we generalize this proof to measure graphs. Immediate corollaries concerning measure preserving actions are given in Section~\ref{section: Applications to measure preserving systems}. We then, in Section~\ref{section: Countable set of translates}, turn to extending the results regarding $c(A,B)$ with $A$ finite, to ones where $A$ is countable. Section~\ref{section: Different summands} is devoted to a proof of a measure-theoretic analogue of Ruzsa's Pl{\"u}nnecke inequality involving different summands. Finally, we prove the correspondence principle for products sets in Section~\ref{section: Applications to Banach Density estimates}.

\textbf{Acknowledgement:} The authors are grateful to Michael Bj{\"o}rklund for many discussions on topics related to the content of the paper and, especially, for the suggestion to work out the case of different summands (Theorem~\ref{thm: Density Plunnecke for different summands}).

\section{Definitions}
\label{section: Definitions}

By a labelled directed graph we mean a tuple $(V,E,A)$ where $V$ and $A$ are sets and $E \subset V \times V \times A$. We regard an element $(v,w,a) \in E$ as an edge directed from $v$ to $w$ and labelled $a$. For subsets $W \subset V$ and labels $a \in A$ the $a$-image and $a$-preimage are defined, respectively, as $$Im_a^{+}(W)=\{ v \in V | (w,v,a) \in E \text{ for some } w \in W  \}$$ and $$Im_a^{-}(W)=\{ v \in V | (v,w,a) \in E \text{ for some } w \in W  \}.$$ That is, the $a$-image of $W$ consists of the vertices that may be approached to from $W$ by walking, in the direction of the orientation, along an edge labelled $a$. Moreover we define for $W\subset V$ the (pre)image $Im^{\pm}(W)=\bigcup_{a \in A} Im_a^{\pm}(W) $. For each integer $h$ we have the $h$-fold image $Im^h(W)$ defined recursively by $Im^0(W)=W$ and $Im^h(W)=Im^{+}(Im^{h-1}(W))$ for $h>0$ and $Im^h(W)=Im^{-}(Im^{h+1}(W))$ for $h<0$. In other words, $Im^h(W)$ consists of all end points of walks with $|h|$ steps that begin at $W$ and agree (resp. disagree) with the orientation of each edge if $h>0$ (resp. $h<0$). Define also the incoming and outgoing degrees of a vetex $v$ as $d^{-}(v)=|\{ (x,y,a) \in E | y=v \}|$ and $d^{+}(v)=|\{ (x,y,a) \in E | x=v \}|$ respectively. Note that $|Im^{\pm}(\{v\})| \leq d^{\pm}(v)$ with strict inequality possible in case of multiple edges between two vertices (of course any two such edges would have different labels). Given an edge $e$ from $v$ to $w$, we will call $v$ the tail, denoted $tail(e)$, and $w$ the head, denoted by $head(e)$. Let $E^+(v)$ denote the edges whose tail is $v$ and $E^-(v)$ those edges whose head is $v$. 

\begin{mydef} A \textit{measure graph} is a tuple $\Gamma = (V,\mathcal{B},\mu,A,E)$ where $(V,\mathcal{B},\mu)$ is a finite measure space (that is, $\mu(V)< \infty$), $A$ is a finite set and $(V,E,A)$ is a labelled directed graph such that

\begin{enumerate}
	\item For each $a \in A$ the sets $$L_a^+(\Gamma) = L^+_a=\{x \in V | \text{ There exists } y \in V \text{ such that } (x,y,a) \in E \}$$ and $$L^{-}_a(\Gamma)=L^{-}_a=\{x \in V | \text{ There exists } y \in V \text{ such that } (y,x,a) \in E \}$$ are measurable.
	\item For $a \in A$ and measurable $W \subset L_a^{\pm }$ we have that $Im_a^{\pm}(W)$ is measurable and $\mu(W)=\mu(Im_a^{\pm}(W))$.
	\item For each label $a \in A$ and vertex $x \in V$ there is at most one outgoing and at most one incoming $a$-labelled edge incident to $x$. That is, $|Im^{\pm}_a(x)|\leq 1$.
\end{enumerate} 

\end{mydef}

\begin{eg}[The $(A,Y,h)$-orbit graph] Given an abelian group $G$ acting on a measure space $(X,\mathcal{B},\mu)$ one may form for each integer $h>0$, finite $A \subset G$ and $Y \in \mathcal{B}$ of finite measure, a measure graph whose underlying vertex set is $\bigsqcup_{k=0}^h A^kY \times \{k\}$ with edge set $$\{ ((x,k),(a.x,k+1),a) | a \in A, x \in A^kY, k=0,1\ldots , h-1 \}.$$ The measure is the restriction of the natural product measure on $X \times \{0,1, \ldots , h\}$.

\end{eg}

Given a labelled graph $\Gamma=(V,E,A)$ and $W \subset V$, the subgraph \textit{induced} by $W$ is the directed labelled graph $(W,E_W,A)$ where $E_W=\{(w_1,w_2,a) | w_1,w_2 \in W, a \in A \text{ and } (w_1,w_2,a) \in E \}$. We say that a subgraph of $\Gamma$ is an \textit{induced subgraph} if it is induced by some subset of $V$.

\begin{eg}\label{eg: induced measure subgraph}  If $\Gamma=(V, \mathcal{B},\mu,A,E)$ is a measure graph and $W \subset V$ is measurable then the subgraph of $\Gamma$ induced by $W$ is a measure graph (with the restricted measure, restricted $\sigma$-algebra and the same edge-label set $A$). Note that the set of vertices with an outgoing edge labelled $a \in A$ is $L_a^{+} \cap W \cap Im_a^{-}(L_a^{-} \cap W)$ and thus is measurable as required.
\end{eg}

Note that the $(A,Y,h)$-orbit graph defined above is a generalization of the commutative addition graph studied in classical Additive Combinatorics, see for example \cite{Ruzsasumsetsandstructure}, \cite{TaoVu}. It is also an example of what is known as a commutative, or Pl{\"u}nnecke, graph which may be defined as follows. 

\begin{mydef} A \textit{layered-graph} is a directed labelled graph $(V,E,A)$ together with a partition $V=V_0 \sqcup V_1 \sqcup \ldots \sqcup V_h$ such that if $e=(x,y,a) \in E$ is a directed edge then $x \in V_i$ and $y \in V_{i+1}$ for some $i \in \{0, \ldots h-1 \}$. We call $V_k$ the $k$-th layer and we say that $(V,E,A)$ is a \textit{$h$-layered} graph (we regard the partition as part of the data of a layered graph). A \textit{semi-commutative} (or \textit{semi-Pl{\"u}nnecke}) graph is a layered graph $(V,E,A)$ such that if $(x,y,a) \in E$ is an edge then there is an injection $\phi:E^{+}(y) \to E^{+}(x)$ such that $(head(\phi(e)),head(e),a) \in E$ for all $e \in E^{+}(y)$. A \textit{commutative}, or \textit{Pl{\"u}nnecke}, graph $\Gamma$ is a directed layered graph such that both $\Gamma$ and the dual graph (the layered graph obtained by reversing edges and the ordering of the layers) $\Gamma^*$ are semi-commutative.

\end{mydef}

\begin{eg} The $(A,Y,h)$-orbit graph defined above is a commutative graph with layering $V=\bigsqcup_{j=0}^h A^jY \times \{j\}$. To check semicommutativity, take a typical edge $((x,j),(ax,j+1),a)$ running from $A^jY \times \{j \}$ to $A^{j+1}Y \times \{j +1\}$ where $0 \leq j \leq h-1$. Then for edges $$e=((ax,j+1),(a'ax,j+2),a') \in E^+((ax,j+1))$$  we may choose $\phi(e)=((x,j),(a'x,j+1),a')$ since, by commutativity of $G$, $a.(a'.x)=(a'a.x)$ and thus there is an $a$-labelled edge from $head(\phi(e))$ to $head(e)$ as required. The semi-commutativity of the dual can be similairly verified.

\end{eg}

The following is an easy exercise in commutative graphs (see \cite{Ruzsasumsetsandstructure}).

\begin{prop} Suppose that $(V,E)$ is a $h$-layered commutative graph with layers $V=V_0 \ldots \sqcup V_h$. Then for $S \subset V_j$ and $T \subset V_k$, where $0 \leq j<k \leq h$, we have that the channel between $S$ and $T$ (that is, the subgraph consisting of all directed paths from $S$ to $T$) is a commutative graph. We denote this subgraph $ch(S,T)$.

\end{prop}

We will be interested in studying channels of an $(A,Y,h)$-orbit graph, it turns out these are measurable.

\begin{lemma} Given an $h$-layered measure graph $\Gamma=(V,\mathcal{B},\mu,A,E)$ with layering\footnote{We always assume implicitly that each layer is measurable.} $V=V_0 \sqcup \ldots \sqcup V_h$ and measurable $S \subset V_i$, $T \subset V_j$ where $0\leq i<j\leq h$ we have that the channel $ch(S,T)$ has measurable vertex set.

\end{lemma}

\textbf{Proof:} Let us denote the vertex set of a subgraph $\Gamma'$ as $V(\Gamma')$. We use induction on $j-i$. The base case $j=i+1$ holds since then $ch(S,T)$ has vertex set $V(ch(S,T))=S \cap Im^- (T) \bigsqcup T \cap Im(S)$. Now suppose that $j-i>1$. Then by the induction hypothesis we have that $ch(S,Im^-(T))$ has measurable vertex set. By the base case, $ch(Im^{-}T,T)$ has measurable vertex set. Now let $U=V(ch(S,Im^- T)) \cap V(ch(Im^- T, T)) \subset V_{j-1}$. Finally we have $V(ch(S,T))=V(ch(S,U)) \cup V(ch(U,T))$ which is measurable again by the induction hypothesis. $\blacksquare$

Note that the previous Lemma and Example~\ref{eg: induced measure subgraph} demonstrate that the channel between two measurable sets may be naturally viewed as a measure graph (as channels are induced subgraphs). 

We now turn to generalizing the notion of the number of edges in a bipartite graph. 

\begin{mydef}  Fix a $1$-layered commutative measure graph $(U,\mathcal{B},\mu,A,E)$ with layering $U=U_0 \sqcup U_1$. Define the \textit{flow} of $\Gamma$ to be the quantity $$ Flow(\Gamma )=\int_{U_0} d^{+}(v) d\mu(v) .$$

\end{mydef}

We now show that the flow behaves nicely and that $d^{\pm}$ is a measurable function. 

\begin{prop}\label{Proposition: Flow} Under the setting of the previous definition, the map $d^+:U \to \mathbb{R}$ is measurable and $Flow( \Gamma )=Flow( \Gamma ^* )$, that is $$ \int_{U_0} d^{+}(v) d\mu(v) = \int_{U_1} d^-(v) d\mu(v).$$

\end{prop}

\textbf{Proof:} Since (by definition of a measure graph) $|Im_a^{\pm}(\{v\})| \leq 1$, we may express $$d^{\pm}=\sum_{a \in A} \chi_{L_a^{\pm}}$$ and thus $d^{\pm}$ is measurable. Consequently we have that $$Flow(\Gamma )=\sum_{a \in A} \mu(L_a^{+})=\sum_{a \in A} \mu (Im_a^+(L_a^{+})) = \sum_{a \in A} \mu (L_a^{-}) = Flow(\Gamma ^*)  $$

as required. $\blacksquare$

\section{Pl{\"u}nnecke's inequality for measure graphs}
\label{section: Plunnecke's inequality for measure graphs}

\begin{mydef} Given a commutative measure graph $\Gamma=(V,\mathcal{B},\mu,A,E)$ with layering $V=V_0 \sqcup \ldots \sqcup V_h$, the magnification ratio of order $j$, where $j \in \{1,\ldots h\}$, is $$ D_j=\inf_{Y \subset V_0, \mu(Y)>0} \frac{\mu(Im^j(Y))}{\mu(Y)}. $$
Moreover, for $C>0$, define the \textit{weight} (corresponding to $C$) to be the measure on $\mathcal{B}$ given by $$w(S)=\sum_{j=0}^h C^{-j}\mu(S \cap V_j) $$ 

for $S \in \mathcal{B}$. Furthermore,  we say that $S \in \mathcal{B}$ is a \textit{cutset} if any path from $V_0$ to $V_h$ intersects $S$ and that $S$ is an \textit{$\epsilon$-minimal cutset} if $S$ is a cutset such that $$w(S) \leq m_0 + \epsilon$$ where $m_0=\inf\{ w(Y) | Y \in \mathcal{B} \text{ is a cutset}  \}.$
\end{mydef}

\begin{lemma}\label{Lemma: U_1 min implies U_0 min} Fix a $2$-layered commutative measure graph $(U,\mathcal{B},\mu,A,E)$ with layering $U=U_0 \sqcup U_1 \sqcup U_2$ and $C>0$. Then if $U_1$ is an $\epsilon$-minimal cutset (with respect to the weight corresponding to $C$), then $U_0$ is an $f(\epsilon)$-minimal cutset where $$f(\epsilon)=\epsilon + 4|A|^2C\epsilon + 4|A|^2\epsilon.$$

\end{lemma}

\textbf{Proof:} Let $m_0=\inf \{ w(S) | S \in \mathcal{B} \text{ is a cutset} \}$. Firstly note that for measurable $S \subset U_1$ we have that $Im(S) \sqcup (U_1 \setminus S)$ is a cutset and thus $m_0 \leq w(Im(S)) + w(U_1 \setminus S)$. On the other hand, since $U_1$ is $\epsilon$-minimal we have that $w(S)+w(U_1 \setminus S) - \epsilon \leq m_0$ and thus $w(S) \leq w(Im(S)) + \epsilon$. A similair argument yields that $w(S) \leq w(Im^{-1}(S)) + \epsilon$. Thus \begin{align} \label{inequality star} C \mu(S) \leq \mu(Im(S)) + C^2 \epsilon  \tag{$\star$} \end{align} and \begin{align}\label{inequality dagger} C^{-1} \mu(S) \leq \mu(Im^{-1}(S)) + \epsilon \tag{$\dagger$} \end{align} for measurable $S \subset U_1$. 

For each integer $i \geq 0$ let  \begin{align*} X_i & =\{ u \in U_1 | d^{-}(v)=i \} \\ Y_i & = \{ u \in U_2 | d^{-}(v)=i \} \\  X'_i & =\{ u \in U_1 | d^{+}(v)=i \}  \\ Y'_i & = \{ u \in U_0 | d^{+}(v)=i \}.  \end{align*}

The $X_i$ are measurable and partition $U_1$. Let $k=|A|$. Define now inductively $T_k=Im(X_k)$ and $T_i=Im(X_{i}) \setminus T_{i+1} $ for $i=k-1, k-2, \ldots ,1$. Note that the $T_i$ partition the set of vertices in $U_2$ that have at least one incoming edge. Moreover, by the definition of a commutativity we have that each vertex in $T_i$ has inwards degree at least $i$ (specifically, this is by the semicommutativity of the dual). Thus we obtain \begin{align} \label{Proof inequality 1} \sum_{i=1}^k i \mu(T_i) \leq  \sum_{i=1}^k Flow(ch(U_1,T_i)) = Flow(ch(U_1,U_2)) \end{align} where the right hand side is well defined since induced subgraphs with measurable vertex sets are measure graphs. From now on we will use the shorthand notation $Flow(U_i,U_j):= Flow(ch(U_i,U_j))$.

Moreover, since $Im(X_j \sqcup \ldots \sqcup X_k)=T_j \sqcup \ldots \sqcup T_k$ we have by (\ref{inequality star}) \begin{align} \label{Proof inequality 2} C \sum_{i=j}^k \mu(X_i) \leq \sum_{i=j}^k \mu(T_i) + C^2 \epsilon . \end{align} Adding these inequalities for $j=1, \ldots k$ we obtain 
\begin{align} \label{Proof inequality 3} C.Flow(U_0,U_1)=C \sum_{i=1}^k i \mu(X_i) \leq \sum_{i=1}^k i \mu(T_i) + kC^2 \epsilon \end{align}

which implies, together with the preceding inequality, that \begin{align} \label{Proof inequality hash} C.Flow(U_0,U_1) \leq Flow(U_1,U_2) + kC^2 \epsilon. \end{align}

We will now apply the same argument to the dual graph to obtain an inequality of the form $$ Flow(U_1,U_2) \leq C. Flow(U_0,U_1) + O_{k,C}(\epsilon).$$ 

To do this, inductively define $T'_k=Im^{-1}(X'_k)$ and $T'_i=Im^{-}(X'_i) \setminus T'_{i+1}$ for $i=k-1,k-2, \ldots ,1$. This time we have, by the duals of the previous arguments (with (\ref{inequality dagger}) in place of (\ref{inequality star})) and the fact that flows are the same for duals (Proposition~\ref{Proposition: Flow}), that \begin{align}\label{Proof inequality 4} \sum_{i=1}^k i \mu(T'_i) \leq Flow(ch_{\Gamma^*}(U_1,U_0))=Flow(U_0,U_1) \end{align}

and, for each $\ell \in \{1,2 \ldots k\}$, 
\begin{align}\label{Proof inequality 5} C^{-1}\sum_{i=\ell}^k \mu(X'_i) \leq \sum_{i=\ell}^k \mu(T'_i) + \epsilon \end{align} 
which, by summing as before, gives $$ C^{-1} Flow(U_1,U_2) \leq Flow(U_0,U_1) + k \epsilon. $$

Thus $Flow(U_1,U_2)$ is close to $C.Flow(U_0,U_1)$. This means that the inequalities above must have been close to being equalities. We will now explicitly estimate how close.  Let us start with the inequality (\ref{Proof inequality 5}). We obtain from it, (\ref{Proof inequality hash}), and (\ref{Proof inequality 4}) that for each $j \in \{1, \ldots k\}$ we have \begin{align*} C^{-1} \sum_{i=j}^k \mu(X'_i) - \sum_{i=j}^k \mu(T'_i) + (k-1) \epsilon & \geq C^{-1} \sum_{i=1}^k i \mu(X'_i) - \sum_{i=1}^k i \mu(T'_i) \\ & \geq C^{-1}.Flow(U_1,U_2) - Flow(U_0,U_1) \\ & \geq - kC\epsilon \end{align*}

where the first inequality is obtained by summing (\ref{Proof inequality 5}) for $\ell \in \{1, 2, \ldots ,k\} \setminus \{j\}$. This finally gives \begin{align} C^{-1}\sum_{i=j}^k \mu(X'_i) - \sum_{i=j}^k \mu(T'_i) \geq -kC \epsilon - (k-1)\epsilon  \end{align}

and so \begin{align} \label{Proof inequality ??} | C^{-1} \sum_{i=j}^k \mu(X'_i) - \sum_{i=j}^k \mu(T'_i) | \leq \max \{kC \epsilon + (k-1)\epsilon, \epsilon  \} \leq kC\epsilon +k\epsilon.  \end{align}

Thus the triangle inequality gives \begin{align} \label{Proof inequality 8} |C^{-1}\mu(X'_i) - \mu(T'_i)| \leq 2kC\epsilon + 2k\epsilon. \end{align}

Now we wish to show that $T'_i$ is approximately $Y'_i$ (that is, they have small symmetric difference). Note that $T'_j \sqcup \ldots \sqcup T'_k \subset Y'_j \sqcup \ldots \sqcup Y'_k$ since, as before, the definition of a commutative graph implies that each vertex in $T_i'$ has outwards degree at least $i$. Thus \begin{align} \sum_{i=j}^k \mu(T'_i) - \sum_{i=j}^k \mu(Y'_i) \leq 0. \end{align}

Combining this with (\ref{Proof inequality 5}) and (\ref{Proof inequality hash}) we have that for each $j \in \{ 1, \ldots k \}$ we have

\begin{align*} \sum_{i=j}^k \mu(T'_i) - \sum_{i=j}^k \mu(Y'_i) & \geq \sum_{i=1}^k i \mu(T'_i) - \sum_{i=1}^k i \mu(Y'_i) \\ & \geq C^{-1}\sum_{i=1}^k i\mu(X'_i) - k\epsilon - \sum_{i=1}^k i\mu(Y'_i) \\ &= C^{-1}.Flow(U_1,U_2) - Flow(U_0,U_1) - k\epsilon \\ & \geq -kC\epsilon - k\epsilon. \end{align*}

Thus $$|\sum_{i=j}^k \mu(T'_i) - \sum_{i=j}^k \mu(Y'_i)| \leq kC\epsilon + k\epsilon$$ from which the trianlge inequality implies \begin{align} |\mu(T'_i) - \mu(Y'_i)| \leq 2kC\epsilon + 2k \epsilon.  \end{align}

Combining this with (\ref{Proof inequality 8}) yields $$|C^{-1} \mu(X'_i) - \mu(Y'_i)| \leq 4kC\epsilon + 4k\epsilon .$$

Finally we get $$|w(U_1)-w(U_0)|=|C^{-1} \sum_{i=1}^k \mu(X'_i) - \sum_{i=1}^k \mu(Y'_i)| \leq 4k^2C\epsilon + 4k^2\epsilon$$

and so in fact $U_0$ is $(\epsilon + 4k^2C\epsilon + 4k^2\epsilon)$-minimal. $\blacksquare$

We will now inductively apply the above lemma to construct $\epsilon$-minimal cutsets that lie in the union of the top and bottom layers. 

\begin{lemma} Suppose that $\Gamma=(V,\mathcal{B},\mu,A,E)$ is a $h$-layered commutative measure graph with layering $X=X_0 \sqcup \ldots \sqcup X_h$. Fix $C>0$ and let $w$ be the weight on $\Gamma$ corresponding to $C$. Then for each $\epsilon>0$ there exists an $\epsilon$-minimal cutset $S \in \mathcal{B}$ such that $S \subset X_0 \sqcup X_h$. 

\end{lemma}

\textbf{Proof:} We will prove, by induction on $j \in \{h-1, \ldots , 1, 0 \}$, that there exists an $\epsilon$-minimal cutset contained in $V_0 \sqcup V_1 \sqcup \ldots \sqcup V_j \sqcup V_h$. The base case $j=h-1$ is clear. Thus fix $\delta>0$ and suppose that $j \in \{1, \ldots , h-1 \}$ and $S \subset V_0 \sqcup V_1 \sqcup \ldots \sqcup V_j \sqcup V_h$ is a $\delta$-minimal cutset. Let $S_i=S \cap X_i$ for $i \in \{1, \ldots ,h\}$. Let $U_0 \subset V_{j-1}$ be those vertices that may be approached to $V_{j-1}$ from $V_0$ along a path that does not intersect $S$. Let $U_2 = V_{j+1} \cap ch(V_{j+1}, V_h \setminus S_h)$ be the set of vertices in $V_{j+1}$ that may be approached to $V_h \setminus S_h$.  We know that $U_2$ is measurable by the measurability of channels and similairly $U_0$ is measurable by an application of the measurability of channels to the subgraph induced by $\bigsqcup_{i=0}^{h} V_i \setminus S_i $. Let $H=ch(U_0,U_2)$ and let $U_1 \subset V_j$ by the vertices in $H$ that lie in $V_j$. Thus $H$ is a $2$-layered measure subgraph of $\Gamma$ that is also commutative. Let us equip $H$ with the measure $C^{-j+1}\mu$ instead of $\mu$, since then the weight function on $H$ corresponding to $C$ agrees with the that of $\Gamma$. 

\underline{Subclaim:} The middle layer $U_1$ is a $\delta$-minimal cutset of $H$.

\begin{figure}[H]
\centering
\includegraphics[scale=0.7]{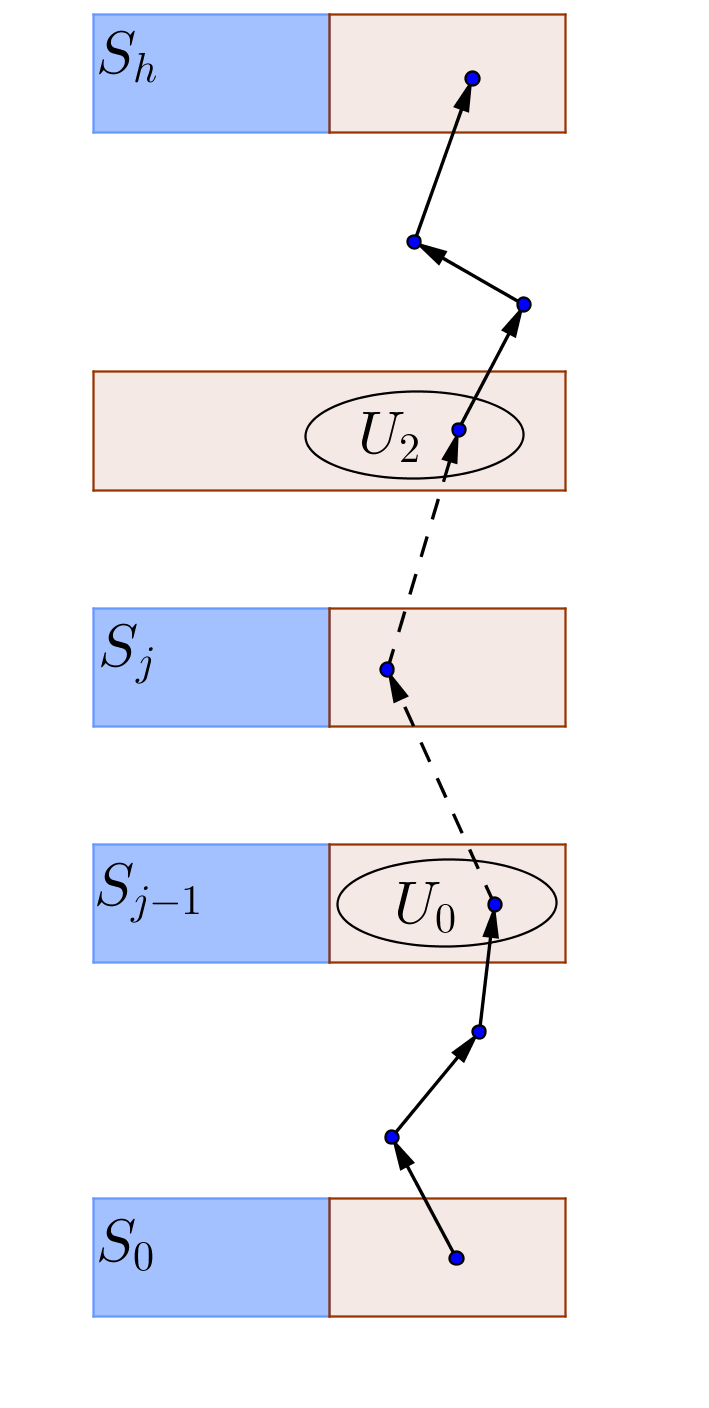}
\caption{A dotted $2$-length path as shown cannot exist as this gives rise to a path from $V_0$ to $V_h$ which avoids $S$, as shown, by the definition of $U_0$ and $U_2$. Thus $U_1 \subset S_j$. One can similairly argue that the $S'$ given in the proof of the subclaim is a cutset.  }
\label{fig:commutative graph subclaim diagram}
\end{figure}

To see this, firstly note that $U_1 \subset S_j$ (see Figure~\ref{fig:commutative graph subclaim diagram}). If $U_1$ is not $\delta$-minimal, then there exists a cutset $T$ in $H$ of weight $w(T) < w(U_1) - \delta \leq w(S_j) - \delta$. But then the set $S'=S_0 \cup \ldots S_{j-1} \cup T \cup S_h$ is a cutset of $\Gamma$ of weight $$w(S') \leq \sum_{i=0}^{j-1} w(S_i) + w(T) + w(S_h) < \sum_{i=0}^{j-1} w(S_i)  + w(S_j) - \delta + w(S_h)=w(S) - \delta,$$ contradicting $S$ being $\delta$-minimal. This proves the subclaim.

Hence we get by Lemma~\ref{Lemma: U_1 min implies U_0 min} that $(S \cup U_0) \setminus S_j$ is a $(\delta + f(\delta))$-minimal cutset, where $f$ is as in the respective lemma (which we may take with the parameters of $\Gamma$, i.e: we consider $H$ as having labelling set $A$ and thus this $f$ does not depend on $H$). Taking $\delta \to 0$ finishes the induction step and hence the proof of this lemma. $\blacksquare$

We are now ready to show that in the case $C=D_h^{1/h}$ the bottom layer is in fact a cutset of minimal weight. 

\begin{corol} Suppose that $\Gamma=(X,\mathcal{B},\mu,A,E)$ is a $h$-layered commutative measure graph with layering $X=X_0 \sqcup \ldots \sqcup X_h$. Suppose that $D_h^{1/h}>0$ and let $w$ be the weight corresponding to $C=D_h^{1/h}$. Then $X_0$ is a cutset of minimal weight. 
\end{corol}

\textbf{Proof:} We want to show that $X_0$ is $\epsilon$-minimal for all $\epsilon>0$. Choose $\epsilon>0$ and by the above lemma an $\epsilon$-minimal cutset $S \subset X_0 \sqcup X_h$. Write $S_i = X_i \cap S$. As $S$ is a cutset we have $Im^h(X_0 \setminus S_0) \subset S_h$ and so $$\mu(S_h) \geq \mu(Im^h(X_0 \setminus S_0)) \geq D_h\mu(X_0 \setminus S_0)=C^h \mu(X_0) - C^h \mu(S_0).$$
Thus $$w(S)=\mu(S_0) + C^{-h}\mu(S_h) \geq \mu(X_0)=w(X_0)$$ and so $w(X_0)$ is $\epsilon$-minimal. $\blacksquare$

We may now finally prove the Pl{\"u}nnecke inequality for measure graphs.

\begin{thm}[Pl{\"u}nnecke inequality for measure graphs] \label{thm: Plunnecke for measure graphs}  Suppose that $\Gamma=(X,\mathcal{B},\mu,A,E)$ is a $h$-layered commutative measure graph with layering $X=X_0 \sqcup \ldots \sqcup X_h$. Then for $j \in \{1, \ldots, h\}$ we have $$D_j^h \geq D_h^j.$$

\end{thm}

\textbf{Proof:} If $D_h=0$ then we are done. If $D_h>0$ then we may set $C=D_h^{1/h}$ and apply the above corollary as follows. For each non-null measurable $Z \subset X_0$ we have that $(X_0 \setminus Z) \sqcup Im^j(Z)$ is a cutset and thus by minimality of $X_0$ we have $$\mu(X_0)=w(X_0) \leq w((X_0 \setminus Z) \sqcup Im^j(Z))=\mu(X_0) - \mu(Z) + D_h^{-j/h}\mu(Im^j(Z))$$ and so $$ D_h^{j/h} \leq \frac{\mu(Im^j(Z))}{\mu(Z)}  $$ which completes the proof as $Z \subset X_0$ was arbitrary. $\blacksquare$

\section{Applications to measure preserving systems} \label{section: Applications to measure preserving systems}

First, we recall the notions of the magnification ratios for the dynamical setting introduced by Bj\"orklund and the second author in \cite{BjorklundFishPlunnecke}.

\begin{mydef} Suppose that $G$ is a countable abelian group acting on a measure space $(X,\mathcal{B},\mu)$. Define for $A \subset G$ and $B \in \mathcal{B}$ of positive finite measure the magnification ratio $$c(A,B) = \inf \{ \frac{\mu(AB')}{\mu(B')} | B' \subset B, \mu(B')>0 \}.$$ Moreover, for $\delta>0$ we may define the $\delta$-heavy magnification ratio $$c_{\delta}(A,B) = \inf \{ \frac{\mu(AB')}{\mu(B')} | B' \subset B, \mu(B') \geq \delta.\mu(B) \}.$$ Furthermore, if $E \subset X$ is measurable then we may define the restricted magnification ratio $$ c(A,B,E) =  \inf \{ \frac{\mu(AB'\setminus E)}{\mu(B')} | B' \subset B, \mu(B')>0 \}. $$ \end{mydef}

By applying the Pl{\"u}nnecke inequality for measure graphs to the case of orbit graphs we obtain the following Pl{\"u}nnecke inequality for measure preserving systems. 

\begin{thm} \label{thm: Plunnecke for G-systems} Suppose that $G$ is a countable abelian group acting on a measure space $(X,\mathcal{B},\mu)$. Then for $A \subset G$ finite and measurable $B \in \mathcal{B}$ of positive finite measure, we have $$c(A^j,B)^{1/j} \geq c(A^k,B)^{1/k}$$ for positive integers $j<k$.

\end{thm}

We may also obtain the $G$-system analogue of a classical restricted addition result. 

\begin{thm} Suppose that $G$ is a countable abelian group acting on a measure space $(X,\mathcal{B},\mu)$. For finite $A \subset G$, measurable $B \subset X$ of positive finite measure and measurable $E \subset X$ we have $$c(A^j,B,A^{j-1}E)^{1/j} \geq c(A^k,B,A^{k-1}E)^{1/k}$$ for positive integer $j<k$.

\end{thm}

\textbf{Proof:} Consider the subgraph of the $(A,B,k)$-orbit graph induced by the subset $$B \times \{0 \} \sqcup \bigsqcup_{j=1}^k (A^jB \setminus A^{j-1}E) \times \{j \}.$$ One may check that this subgraph is indeed commutative (see \cite{Ruzsasumsetsandstructure}). $\blacksquare$

\section{Countable set of translates}
\label{section: Countable set of translates}

The inequalities established in Section~\ref{section: Applications to measure preserving systems} required the set of translates $A \subset G$ to be finite. We now turn to extending Theorem~\ref{thm: Plunnecke for G-systems} to the case where $A$ is countable. We use the techniques developed by Bj\"orklund and the second author in \cite{BjorklundFishPlunnecke}.

The following proposition is analogous to Proposition 2.2 in \cite{BjorklundFishPlunnecke}. 

\begin{prop} Suppose that $G$ is an abelian group acting on a probability space $(X,\mathcal{B},\mu)$ and fix a finite $A \subset G$ and non-null $B \in \mathcal{B}$ together with a $0 < \delta <1 $ and positive integers $j \leq k$. If $B' \subset B$ is measurable and satisfies \begin{align} \label{Heavy hypothesis}\left( \frac{\mu(A^kB')}{\mu(B')} \right) ^{1/k} \leq (1-\delta)^{-1/j} \left( \frac{\mu(A^jB)}{\mu(B)}  \right) ^{1/j}. \end{align}

Then $\mu(B') \geq \delta. \mu(B)$ or there exists $B' \subset B'' \subset B$ such that $\mu(B''\setminus B)>0$ and $B''$ satisfies (\ref{Heavy hypothesis}).

\end{prop}

\textbf{Proof:} Firstly we note that if the hypothesis holds for $B'=B_1$ and $B'=B_2$ with $B_1$ and $B_2$ disjoint, then it holds for $B_1 \sqcup B_2$ since the hypothesis may be rewritten as the inequality $$\mu(A^kB') \leq (1-\delta)^{-k/j} \left( \frac{\mu(A^jB)}{\mu(B)}  \right) ^{k/j} \mu(B').$$

By Theorem~\ref{thm: Plunnecke for G-systems} we know that there exists non-null measurable $B' \subset B$ such that (\ref{Heavy hypothesis}) is satisfied. Suppose that $\mu(B') < \delta . \mu(B)$, thus we wish to construct a strictly larger $B'' \supset B$ that satisfies (\ref{Heavy hypothesis}) and is contained in $B$. Set $B_0 = B \setminus B'$. We have that $$ \frac{\mu(B_0)}{\mu(B)} (1-\delta)^{-1} > 1 $$ and thus there exists $B_0' \subset B_0$ such that \begin{align*} \left( \frac{\mu(A^kB'_0)}{\mu(B'_0)} \right) ^{1/k} & \leq \left( \frac{\mu(B_0)}{\mu(B)}  (1-\delta)^{-1} \right)^{1/j} \left( \frac{\mu(A^jB_0)}{\mu(B_0)}  \right) ^{1/j} \\ & = (1-\delta)^{-1/j} \left( \frac{\mu(A^jB_0)}{\mu(B)}  \right) ^{1/j} \\ & \leq (1-\delta)^{-1/j} \left( \frac{\mu(A^jB)}{\mu(B)}  \right) ^{1/j} \end{align*}

and thus we may set $B''=B' \sqcup B_0'$. $\blacksquare$

We will now apply the above lemma to construct a set $B' \subset B$ such that $\mu(B') \geq \delta.\mu(B)$ and (\ref{Heavy hypothesis}) holds. The idea is to choose a set $B' \subset B$ that satisfies (\ref{Heavy hypothesis}) and that is maximal in the sense that $B$ does not contain any measurable $B'' \supset B'$ of strictly larger measure that satisfies (\ref{Heavy hypothesis}). Such a set would have to necessarily satisfy $\mu(B') \geq \delta. \mu(B)$. The existence of such a maximal $B'$ follows from the continuity of measure together with the following easy lemma on monotone classes. 

\begin{lemma}\label{lemma: monotone maximality} Suppose that $(X,\mathcal{B}, \mu)$ is a finite measure space and $\mathcal{M} \subset \mathcal{B}$ is non-empty and closed under countable nested unions (that is, if $M_i \in \mathcal{M}$ with $M_i \subset M_{i+1}$ then $\bigcup_{i=1}^{\infty} M_i \in \mathcal{M}$). Then there exists $M \in \mathcal{M}$ such that $\mu(M)=\mu(M')$ for all $M' \in \mathcal{M}$ with $M \subset M'$.

\end{lemma}

\textbf{Proof:} For $M \in \mathcal{M}$ let $s(M)=\sup \{\mu(M') | M' \in \mathcal{M}, M \subset M' \}$. Choose $M_1 \in \mathcal{M}$. Now inductively choose $M_{n+1} \in \mathcal{M}$ such that $M_n \subset M_{n+1}$ and $\mu(M_{n+1}) \geq \frac{\mu(M_n) + s(M_n)}{2}$. Let $M= \bigcup_{n=1}^{\infty}M_n$. We claim that $\mu(M)=s(M)$. To see this, note that $\mu(M_n) \to \mu(M)$ and $s(M_n) \geq s(M)$. Thus $$ \mu(M) = \lim_{n \to \infty} \mu(M_{n+1}) \geq \limsup_{n\to \infty} \frac{\mu(M_n) + s(M_n)}{2} \geq \limsup_{n\to \infty} \frac{\mu(M_n) + s(M)}{2} = \frac{ \mu (M) + s(M) }{2}$$ and thus $\mu(M) \geq s(M)$ as required. $\blacksquare$

If we set $\mathcal{M} = \{ B' \subset B | B' \text{ satisfies (\ref{Heavy hypothesis})} \}$ then we see that $\mathcal{M}$ is non-empty by Theorem~\ref{thm: Plunnecke for G-systems} and is closed under countable nested unions by the continuity of measure. Thus by the discussion above we obtain a $B' \subset B$ such that $\mu(B') \geq \delta.\mu(B)$ and (\ref{Heavy hypothesis}) holds. Consequently we have shown 

\begin{lemma} Suppose that $G$ is an abelian group acting on a probability space $(X,\mathcal{B},\mu)$ and fix a finite $A \subset G$ and non-null $B \in \mathcal{B}$ together with a $0 < \delta <1 $ and positive integers $j \leq k$. Then $$c_{\delta}(A^k,B)^{1/k} \leq (1-\delta)^{-1/j} \left( \frac{\mu(A^jB)}{\mu(B)}  \right) ^{1/j}.$$

\end{lemma}

We may now obtain our first result about the case where $A \subset G$ is not necessarily finite. 

\begin{lemma} \label{lemma: sup c_d upper bound} Suppose that $G$ is an abelian group acting on a probability space $(X,\mathcal{B},\mu)$ and fix a (not necessarily finite) set $A \subset G$ and  non-null $B \in \mathcal{B}$ together with a $0 < \delta <1 $ and positive integers $j \leq k$. Then 

$$ \sup \{ c_{\delta}(A',B)^{1/k} | A' \subset A^k, \text{ } A' \text{ is finite} \} \leq (1-\delta)^{-1/j} \left( \frac{\mu(A^jB)}{\mu(B)}  \right) ^{1/j}. $$

\end{lemma}

\textbf{Proof:} If $A' \subset A^k$ is finite then one may choose a finite $A_0 \subset A$ such that $A' \subset A_0^k$. Consequently $$c_{\delta}(A',B)^{1/k} \leq c_{\delta}(A_0^k,B)^{1/k} \leq (1-\delta)^{-1/j} \left( \frac{\mu(A_0^jB)}{\mu(B)}  \right) ^{1/j} \leq (1-\delta)^{-1/j} \left( \frac{\mu(A^jB)}{\mu(B)}  \right) ^{1/j} $$ and so as $A'$ was arbitrary this completes the proof. $\blacksquare$

The next non-trivial result due to Bj\"orklund and Fish allows us to extend the Pl{\"u}nnecke inequalities for a finite set of translates (Theorem~\ref{thm: Plunnecke for G-systems}) to the case of an infinite set of translates.   

\begin{thm}[Proposition 4.1 of \cite{BjorklundFishPlunnecke}]  Suppose that $G$ is a countable group acting on a probability space $(X,\mathcal{B},\mu)$ such that $L^2(X,\mathcal{B},\mu)$ is separable and fix a (not necessarily finite) set $A \subset G$ and  non-null $B \in \mathcal{B}$ together with a $0 < \delta <1 $. Then $$c(A,B) \leq \sup \{ c_{\delta}(A',B) | A' \subset A, \text{ } A' \text{ is finite} \}.$$

\end{thm} 

\begin{thm}(Pl\"unnecke inequalities for an infinite set of translates) \label{thm: Plunnecke for infinite translates} Suppose that $G$ is a countable abelian group acting on a probability space $(X,\mathcal{B},\mu)$ such that $L^2(X,\mathcal{B},\mu)$ is separable and fix a (not necessarily finite) set $A \subset G$ and non-null $B \in \mathcal{B}$ together with a $0 < \delta <1 $ and positive integers $j \leq k$. Then 

$$ c(A^k,B)^{1/k} \leq c(A^j,B)^{1/j} \leq \left( \frac{\mu(A^jB)}{\mu(B)}  \right) ^{1/j}. $$

\end{thm} 

\textbf{Proof:} By the previous two results we obtain for each $\delta>0$ the inequalities $$c(A^k,B)^{1/k} \leq \sup \{ c_{\delta}(A',B)^{1/k} | A' \subset A^k, \text{ } A' \text{ is finite} \} \leq (1-\delta)^{-1/j} \left( \frac{\mu(A^jB)}{\mu(B)}  \right) ^{1/j}.  $$

Taking $\delta \to 0$ gives $$c(A^k,B)^{1/k} \leq \left( \frac{\mu(A^jB)}{\mu(B)}  \right) ^{1/j}. $$

Now applying this to non-null $B_i \subset B$ such that $$\frac{\mu(A^jB_i)}{\mu(B_i)} \to c(A^j,B)$$ gives  $$ c(A^k,B)^{1/k} \leq c(A^k,B_i)^{1/k} \leq \left( \frac{\mu(A^jB_i)}{\mu(B_i)}  \right) ^{1/j} \to c(A^j, B)^{1/j} $$ as desired. $\blacksquare$

\section{Different summands}
\label{section: Different summands}

Given measure preserving actions $G \curvearrowright (X, \mathcal{B},\mu)$ and $G' \curvearrowright (X', \mathcal{B}',\mu')$ one can form the measure preserving product action $G \oplus G \curvearrowright (X \times X, \mathcal{B} \otimes \mathcal{B}', \mu \times \mu')$ given by $(g,g').(x,x')=(g.x,g'.x')$. We will now verify that the corresponding multiplication ratios are multiplicative.

\begin{lemma} \label{Lemma: multiplicativity of mag ratios} Suppose that $G$ and $G'$ are countable groups acting on probability spaces  $(X, \mathcal{B},\mu)$ and  $(X', \mathcal{B}',\mu')$ respectively. Then for $A \subset G$, $A' \subset G'$ and non-null $B \subset X$, $B' \subset X'$ we have $$c(A,B) c(A',B')=c(A \times A', B \times B'). $$

\end{lemma}

\textbf{Proof:} For non-null $B_0 \subset B$ and $B'_0 \subset B'$ we have \begin{align*} c(A\times A', B \times B') & \leq \frac{\mu \times \mu' (A \times A'.B_0 \times B'_0)}{\mu \times\mu' (B_0 \times B'_0)} \\ & = \frac{\mu(A.B_0)}{\mu(B_0)} \frac{\mu(A'.B'_0)}{\mu(B'_0)}. \end{align*} By selecting appropriate $B_0$ and $B'_0$, the right hand side may be made arbitrarily close to $c(A,B) c(A',B')$ and thus $$ c(A,B) c(A',B') \geq c(A \times A', B \times B'). $$ 

We now aim to show the reverse inequality. For $U \subset X \times X'$ and $(x_0,x'_0) \in X \times  X'$ let $$U_{x_0}= \{ x' \in X' | (x_0,x') \in U \}$$ and $$U^{x'_0} = \{ x \in X | (x,x_0') \in U \}.$$ Also, let $\mathcal{B}(V)$ denote the measurable subsets of $V$ where $V$ is any subset of a measurable space. Define $$\phi : \mathcal{B} ( X \times B') \to \mathcal{B} (X \times X')$$ by $$\phi(U) = \bigsqcup_{x \in X} \{x \} \times A' U_x = (\{1_G\} \times A'). U . $$

By Fubini's theorem we have 

\begin{align*} \mu \times \mu' ( \phi(U) ) & = \int_X \mu'(A'U_x) d\mu(x) \\ & \geq \int_X c(A',B') \mu'(U_x) d\mu(x) \\ & = c(A',B') \int_X \mu' (U_x)d\mu(x)  \\ & = c(A',B') \cdot \mu \times \mu' (U).  \end{align*}

Thus $\mu \times \mu' ( \phi (U )) \geq c(A',B') \cdot \mu \times \mu' (U)$ for $U \subset X \times B'$. We may reverse the role of co-ordinates to obtain a similair inequality, from which we finally get that \begin{align*} \mu \times \mu' ( (A \times A') .U ) & = \mu \times \mu' ((\{1_G \} \times A')(A \times \{1_{G'} \}) . U)) \\ & \geq c(A',B') \cdot \mu \times \mu'  ((A \times \{1_{G'} \}) . U)) \\ & \geq c(A',B') c(A,B) \cdot \mu \times \mu' (U) \end{align*} for $U \subset B \times B'$. This implies that  $c(A,B) c(A',B') \leq c(A \times A', B \times B')$, as required. $\blacksquare$

\begin{prop}\label{prop: different summands} Suppose that $G$ is a countable abelian group acting on a probability space $(X, \mathcal{B}, \mu)$. Then for $A_1, A_2, \ldots A_k \subset G$ and non-null $B \in \mathcal{B}$ we have $$ c(A_1 \ldots A_k, B) \leq \prod_{i=1}^k \frac{\mu(A_i B)}{\mu(B)}.$$

\end{prop}

\textbf{Proof:} Choose rational numbers $$\alpha_i > \frac{\mu(A_iB)}{\mu(B)}$$ and choose $n \in \mathbb{Z}_{>0}$ such that for each $i \in \{1,\ldots k\}$ we have $$n_i := \frac{n}{\alpha_i} \in \mathbb{Z}_{>0}.$$ Suppose that there exists $T_i \subset G$ with $|T_i|=n_i$ such that the map \begin{align*} T_1 \times \ldots \times T_k \times (A_1A_2\ldots A_kB) & \to X \\ (t_1, \ldots ,t_k, y) & \mapsto t_1\ldots t_k y \end{align*} is injective. We may assume, without loss of generality, that we are in this case by naturally embedding $G \hookrightarrow G \oplus \mathbb{Z}/N\mathbb{Z}$ and $X \hookrightarrow X \times \mathbb{Z}/N\mathbb{Z}$ and replacing the measure preserving system $G \curvearrowright X$ with the product measure preserving system $G \oplus \mathbb{Z}/N\mathbb{Z} \curvearrowright X \times \mathbb{Z}/N\mathbb{Z}$, for large enough $N$. Let $A=\bigcup_{i=1}^k A_iT_i$ and notice that  $$ \mu(AB) \leq \sum_{i=1}^k \mu(A_iT_iB) \leq \sum_{i=1}^k n_i \cdot \mu(A_iB) < \mu(B) \sum_{i=1}^k n_i \alpha_i = k \cdot n \cdot   \mu(B)$$

and thus, by Theorem~\ref{thm: Plunnecke for infinite translates}, we obtain non-null $B' \subset B$ such that $$ \mu(A^kB') \leq \left( k \cdot n \right)^k \mu(B').  $$ 

However, by the injection above, we have $$\mu(A^kB') \geq \mu( T_1 \ldots T_k A_1A_2 \ldots A_k B') = \left( \prod_{i=1}^k n_i \right) \mu(A_1 \ldots A_k B').$$ Combining the previous two inequalities gives 
$$c(A_1\ldots A_k,B) \leq \frac{\mu(A_1 \ldots A_k B')}{\mu(B')} \leq k^k \prod_{i=1}^k \alpha_i  .$$

Since the $\alpha_i > \frac{\mu(A_iB)}{\mu(B)}$ were arbitrary rational numbers, we obtain 

\begin{align}\label{k^k bound} c(A_1\ldots A_k,B) \leq k^k \prod_{i=1}^k \frac{\mu(A_iB)}{\mu(B)} . \end{align}

We now wish to remove the $k^k$ constant. This may be done by consdering a large cartesian power, as follows. Let us denote $V^{\times m}=V \times V \ldots \times V$ and $V^{\oplus m}=V\oplus V\ldots \oplus V$, where $m$ factors are present. For each positive integer $m$, an application of (\ref{k^k bound}) and Lemma~\ref{Lemma: multiplicativity of mag ratios} to the sets $A_1^{\oplus m}\ldots A_k^{\oplus m} = (A_1A_2\ldots A_k) ^{\oplus m} \subset G^{\oplus m}$ and $B^{\times m} \subset X^{\times m}$ gives $$ c(A_1\ldots A_k, B) \leq k^{k/m}  \prod_{i=1}^k \frac{\mu(A_iB)}{\mu(B)}. $$

Taking the limit $m \to \infty$ gives the desired result. $\blacksquare$

\section{Correspondence principle for product sets}
\label{section: Applications to Banach Density estimates}

We will now establish a Furstenberg correspondence principle for product sets. The first appearance of the correspondence principle for product sets was in \cite{countableamenablegroups2010}. The principle appearing in this paper is due to Bj\"orklund and Fish, and appeared in \cite{BjorklundFishProductSet}. Given a countable group $G$ and $B \subset G$, we define the Furstenberg \textit{$G$-system corresponding to $B$} to be the topological $G$-system $G \curvearrowright X$, i.e., G acts on compact metric space X by homeomorphisms, given by the following construction. Let $X_0=\{0,1\}^G$ be the space of all sequences indexed by $G$ equipped with the product topology. Let $z \in \{0,1\}^G$ be the indicator function of $B$, that is $z_g=1$ if and only if $g \in B$. Note that there is a natural action of $G$ on $X_0$ given by $$ (g.x)_h=x_{gh}$$ for $g,h\in G$ and $(x_k)_{k \in G} \in X_0$.

Let $X = \overline{Gz}=\overline{ \{ gz \text{ } | \text{ } g \in G \} }$ be the closure of the orbit of $z$. Note that $X$ is $G$-invariant. This defines the system corresponding to $B$. Moreover, we define the \textit{clopen set corresponding to $B$} to be the set $$\widetilde{B}=\{ x \in X \text{ } |\text{ }x_1 = 1\},$$ which is a clopen subset of $X$.

\begin{lemma} Suppose that $G$ is a countable abelian group, $B \subset G$ and $X$ is the $G$-system corresponding to $B$. Suppose that $\mu$ is a $G$-invariant Borel probability measure on $X$. Then for finite $A_0 \subset G$ we have that $$ d^*(A_0B) \geq \mu(A_0\widetilde{B}) \geq d_*(A_0B). $$ 

\end{lemma}

\textbf{Proof:} Note that it suffices to prove this for $\mu$ ergodic by either of the following arguments. Suppose that the result holds for all ergodic $\mu$, and thus holds for all convex combinations of ergodic measures. It is a well known fact that the extreme $G$-invariant measures are precisely the ergodic measures. The Krein-Milman theorem therefore implies that any $G$-invariant probability measure is in the weak$^*$ closure of the set of all convex combinations of ergodic measures. Since the map $\nu \mapsto \nu(A_0\widetilde{B})$ is weak$^*$ continuous (since $A_0\widetilde{B}$ is clopen), we obtain the result. Alternatively, one may use Bauer's maximum principle (instead of the Krein-Milman) which says that the maximum (resp. minimum) of the the map $\nu \mapsto \nu (A_0\widetilde{B})$ is attained at an extremal (hence ergodic) measure, say $\mu^*$ (resp. $\mu_*$), and thus for any $G$-invariant $\mu$ we have $$d_* (A_0B)\leq \mu_*(A_0\widetilde{B}) \leq \mu(A_0\widetilde{B}) \leq \mu^*(A_0\widetilde{B}) \leq d^*(A_0B).$$ Now we turn to the proof of the Lemma under the assumption that $\mu$ is ergodic. Given any F{\o}lner sequence $(F_n)_{n \in \mathbb{N}}$ and continuous $f \in C(X)$ we have, by the Von Neumann mean ergodic theorem, that \begin{align}\label{mean limit}\frac{1}{|F_N|}\sum_{g \in F_N} f \circ g \to \int f d\mu \end{align} in the $L^2$-norm. We may then, by the Borel-Cantelli lemma, pass to a subsequence of $(F_n)_{n \in \mathbb{N}}$ to obtain almost everywhere pointwise convergence in (\ref{mean limit}). In particular we may apply this result to $f=\chi_{A_0\widetilde{B}}$ and get a F{\o}lner sequence $(F_n)_{n \in \mathbb{N}}$ such that $$\frac{1}{|F_N|}\sum_{g \in F_N} \chi_{A_0\widetilde{B}}(gx) \to \mu(A_0\widetilde{B})$$ for some $x \in X$. Now fix $N \in \mathbb{N}$ and note that since $X=\overline{Gz}$ we have $h_iz \to x$ for some $h_i \in G$ and thus $\chi_{A_0\widetilde{B}}(gh_iz) \to \chi_{A_0\widetilde{B}}(gx)$ for each $g \in F_N$. Therefore, for some large $M$, we have for $q_N:=h_M$ the equality $$\frac{1}{|F_N|}\sum_{g \in F_N} \chi_{A_0\widetilde{B}}(gx) = \frac{1}{|F_N|}\sum_{g \in F_N} \chi_{A_0\widetilde{B}}(gq_Nz) = \frac{|q_N F_N \cap A_0B|}{|q_NF_N|}  $$

where the second equality is obtained from the fact that, by construction of the corresponding system and clopen set, we have $gq_Nz \in A_0\widetilde{B}$ if and only if $gq_N \in A_0B$. Since $(q_NF_N)_{N \in \mathbb{N}}$ is a F{\o}lner sequence, the limit (as $N \to \infty$) of this quantity must be between $d_* (A_0B)$ and $d^*(A_0B)$. $\blacksquare$ 

In fact, notice that the inequality $d^*(A_0B) \geq \mu(A_0\widetilde{B})$ in the previous lemma is also true for infinite $A_0$ since we can always write $A_0$ as an increasing union $A_1 \subset A_2 \ldots $ of finite sets and thus $$d^*(A_0B)\geq d^*(A_kB) \geq \mu(A_k\widetilde{B}) \to \mu(A\widetilde{B})\text{ as } k \to \infty.$$

\begin{prop}[Correspondence principle for product sets \cite{BjorklundFishPlunnecke}] Suppose that $G$ is a countable abelian group and $A,B \subset G$. Then there exists a compact metrizable space $X$ on which $G$ acts by homeomorphisms such that there exist $G$-invariant ergodic Borel probability measures $\mu, \nu$ on $X$ together with a clopen $\widetilde{B} \subset X$ such that \begin{align*}  d^* (B) &= \mu(\widetilde{B}) \\ d^*(AB) &\geq \mu(A\widetilde{B}) \\ d_* (B)  &\leq \nu(\widetilde{B}) \\ d_*(AB)  &\geq \nu(A\widetilde{B}). \end{align*} 

\end{prop}

\textbf{Proof:} The space $X$ and clopen set $\widetilde{B}$ will be those coming from the correspondence. Note that the second and third inequalities are satisfied for all $\mu,\nu$ by the lemma above. Moreover, this lemma shows that  $d^* (B) \geq \mu(\widetilde{B})$, for all $\mu$. Therefore to construct $\mu$ satisfying the first equality, it is enough to construct a not necessarily ergodic $\mu$ and then apply Bauer's maximum principle. Let $z$ be as in the construction of the correspondence. Choose a F{\o}lner sequence $F_n \subset G$ such that $$\frac{|F_n \cap B|}{|F_n|} \rightarrow d^*(B).$$ Consider now the following averages of point mass measures $$ \mu_n = \frac{1}{|F_n|} \sum_{g \in F_n} \delta_{g.z}$$ and let $\mu = \lim_{k \to \infty} \mu_{n_k}$ be a weak$^*$ limit of a subsequence of these. Since $(F_n)_{n \in \mathbb{N}}$ is F{\o}lner, we have the $\mu$ is $G$-invariant. Note that for $C \subset G$ $$\mu_n(C\widetilde{B}) = \frac{|F_n \cap CB |}{|F_n|}. $$ In particular, the $C=\{1 \}$ case shows that the choice of F{\o}lner sequence, together with the fact that $\widetilde{B}$ is clopen, implies that $d^{*}(B)=\mu(\widetilde{B})$. Now we turn to dealing with the final inequality. To construct such a $\nu$, it is enough to construct such a not necessarily ergodic $\nu$ by the following argument. The map $\nu \mapsto \nu(A\widetilde{B})$ is weak$^*$ lower semicontinuous since $A\widetilde{B}$ is open. Thus Bauer's minimum (but not maximum) principle applies and thus if at least one not necessarily ergodic $\nu$ satisfies the final inequality, then some ergodic minimizer does too. Write $A$ as a union of an increasing sequence of finite sets $A_1 \subset A_2 \subset \ldots $ and choose a F{\o}lner sequence $E_m \subset G$ such that $$\frac{|E_m \cap AB|}{|E_m|} \to d_{*}(AB) \text{ as } m \to \infty.$$

As before, we have that the averages $$ \nu_m = \frac{1}{|E_m|} \sum_{g \in E_n} \delta_{g.z} $$ have a weak$^*$ convergent subsequence $\nu_{m_j} \to \nu$. Since each $A_k\widetilde{B}$ is clopen, we have $$d_*(AB) = \lim_{j \to \infty} \frac{|E_{m_j} \cap AB|}{|E_{m_j}|} \geq \lim_{j \to \infty} \frac{|E_{m_j} \cap A_kB|}{|E_{m_j}|} = \lim_{j \to \infty} \nu_{m_j} (A_k\widetilde{B}) = \nu(A_kB) \to \nu(A\widetilde{B}) \text{ as } k \to \infty.  $$
$\blacksquare$

The following statement was proven in \cite{BjorklundFishPlunnecke}.

\begin{lemma}[\cite{BjorklundFishPlunnecke}]\label{lemma: Standard density estimate from ergodic action} Suppose that $G$ is a countable abelian group that acts ergodically on a probability space $(X,\mathcal{B},\mu)$ such that $L^2=L^2(X,\mathcal{B},\mu)$ is separable (for instance, a Borel probability space). Then $$d^*(A) \leq \mu(AB) $$ for $A \subset G$ and non-null $B \in \mathcal{B}$.

\end{lemma}

\textbf{Proof:} Choose a F{\o}lner sequence $(F_n)_{n \in \mathbb{N}}$ such that $$ \lim_{n \to \infty} \frac{|F_n \cap A|}{|F_n|} = d^*(A).$$ Define $$f_N = \frac{1}{|F_N|} \sum_{g \in F_N} \chi_{g^{-1}AB} \in L^2.$$

As $\| f_N \|_2 \leq 1$ and $L^2$ is separable (and thus has unit ball compact metrizable in the weak topology), we may pass to a subsequence of $(F_n)_{n \in \mathbb{N}}$ such that $f_N$ converges weakly to some $f \in L^2$. But $f$ is $G$-invariant and thus constant by ergodicity. Therefore $$\mu(AB)=\langle f_N,1 \rangle \to \langle f,1 \rangle=f$$ and so in fact $f$ is the constant function $\mu(AB)$.

Notice that for $b \in B$ we have \begin{align*} f_N(b) & = \frac{1}{|F_N|} | \{ g \in F_N | b \in g^{-1}AB \}|
\\ & = \frac{1}{|F_N|} | \{ g \in F_N | gb \in AB \}| \\ & \geq \frac{1}{|F_N|} | \{ g \in F_N | g \in A \}| \end{align*}

and so $d^*(A) \leq \liminf_{N \to \infty} f_N(b).$ In other words, we have $$\chi_B d^*(A) \leq \chi_B \liminf_{N \to \infty} f_N. $$ Integrating this inequality and applying Fatou's lemma yields \begin{align*}\mu(B) d^*(A) & \leq \int \chi_B \liminf_{N \to \infty} f_N \\ & \leq \liminf_{N \to \infty} \int \chi_B f_N \\ & = \langle \chi_B, f \rangle \\ & = \mu(B)\mu(AB). \end{align*}

As $\mu(B)>0$, we obtain the desired inequality. $\blacksquare$

\bibliographystyle{amsplain}
\bibliography{mybib}

\end{document}